\documentclass{article}
\usepackage{amssymb,amsmath,amsfonts,amsthm}
\usepackage{latexsym}
\usepackage{mathrsfs}
\usepackage{graphics}
\usepackage{tikz}
\usetikzlibrary{shapes,arrows,calc}
\usepackage{indentfirst}

\usepackage{hyperref}
\usepackage{enumerate}
\usepackage{mathtools}
\usepackage{comment}
\usepackage{xparse}
\usepackage{bm}
\usepackage{subcaption}
\usepackage{caption}
\usepackage{graphicx} 
\usetikzlibrary{positioning} 
\hypersetup{colorlinks = true, linkcolor = blue, citecolor = blue, urlcolor = blue}
\usepackage{float}
\usepackage[section]{placeins}
\usepackage{lineno}

\allowdisplaybreaks

\newtheorem*{thm*}{Theorem}
\newtheorem{thm}{Theorem}
\newtheorem{lem}[thm]{Lemma}
\newtheorem{pro}[thm]{Proposition}
\newtheorem{obs}[thm]{Observation}
\newtheorem{cor}[thm]{Corollary}
\newtheorem{conj}[thm]{Conjecture}

\newcommand{\N}{\mathbb{N}}

\newcommand{\col}{\mathrm{col}}

\usepackage[margin=1in]{geometry}
\usepackage{amsmath,amssymb}
\usepackage{amsthm}

\newtheorem*{claim*}{Claim}

\usepackage{xparse} 
\usepackage{etoolbox}
\NewDocumentCommand{\angles}{mmO{1}}{%
  \ensuremath{%
    \ifcase\numexpr#3\relax
      \PackageError{angles}{Nesting level must be >= 1}{You passed #3}%
    \or
      #1\langle #2 \rangle
    \or
      #1\langle\langle #2 \rangle\rangle
    \or
      #1\langle\langle\langle #2 \rangle\rangle\rangle
    \else
      \PackageError{angles}{Unsupported nesting level #3}{Only nesting levels 1--3 are supported}
    \fi
  }%
}
\usetikzlibrary{fit,matrix,positioning,backgrounds,shapes.geometric}
\newsavebox{\KTwoFourFigureBox}

\begin{document}
\title{An Upper Bound on the Proportional Choice Number of a Graph}

\author{Kevin Grace \thanks{Department of Mathematics and Statistics, University of South Alabama, Mobile, AL, USA (kevingrace@southalabama.edu)} \and Jeffrey A. Mudrock \thanks{Department of Mathematics and Statistics, University of South Alabama, Mobile, AL, USA (mudrock@southalabama.edu)}
}

\maketitle
\begin{abstract}

Proportional choosability is a list analogue of equitable coloring introduced in 2019 that imposes both upper and lower bounds on color class sizes based on the availability of each color. The proportional choice number of a graph $G$, denoted by $\chi_{pc}(G)$, is the smallest positive integer $k$ such that $G$ is proportionally $k$-choosable. A recently proposed conjecture states that for any graph $G$ with at least one edge, $\chi_{pc}(G) \leq \Delta(G) + \col(G) - 1$, where $\Delta(G)$ is the maximum degree of $G$ and $\col(G)$ is its coloring number. In this paper, we prove this conjecture using a generalization of Edmonds' matroid intersection theorem due to Aharoni and Berger, in which one of the matroids is replaced with an abstract simplicial complex. As a consequence, we can completely determine the proportional choice number of every path: $\chi_{pc}(P_1)=1$, $\chi_{pc}(P_n)=2$ for $2 \leq n \leq 5$, and $\chi_{pc}(P_n)=3$ for $n \geq 6$.  This resolves a question about paths that has remained open since the introduction of proportional choosability.



\end{abstract}

\section{Introduction} \label{intro}

In this paper, all graphs are finite and simple, and unless otherwise stated, they are nonempty. We allow the empty graph when it arises from deleting vertices. Generally speaking, we follow West~\cite{W01} for terminology and notation. The set of natural numbers is $\N = \{1,2,3,\ldots\}$. For $m \in \N$, we write $[m]$ for the set $\{1,\ldots,m\}$, and we let $[0]=\emptyset$. We write $K_{l,n}$ for complete bipartite graphs with partite sets of size $l$ and $n$, and we write $P_n$ for the path on $n$ vertices. When we write $G=qK_{l,n}$, we mean that $G$ is the disjoint union of $q$ copies of $K_{l,n}$.  For a graph $G$, we use $V(G)$ and $E(G)$ for its vertex set and edge set, respectively. We use $d_G(v)$ for the degree of a vertex $v$ in $G$, and we use $\Delta(G)$ for the maximum degree of the vertices in $G$. When $u,v \in V(G)$ are adjacent in $G$, we write $uv$ or $vu$ for the edge connecting $u$ and $v$. For each $v \in V(G)$, we let $N_G(v)$ denote the set of neighbors of $v$ in $G$, and we let $N_G[v]$ denote the closed neighborhood of $v$ in $G$.

Suppose $G$ is a graph and $X \subseteq V(G)$. We let $G[X]$ denote the subgraph of $G$ induced by $X$, and we let $G-X$ denote the graph obtained from $G$ by deleting all vertices in $X$ and their incident edges. If $e \in E(G)$, we let $G-e$ denote the graph obtained from $G$ by deleting the edge $e$.  Finally, the \emph{coloring number} of a nonempty graph $G$, denoted by $\col(G)$, is the smallest positive integer $k$ such that every nonempty subgraph of $G$ has a vertex of degree at most $k-1$.

\subsection{Equitable List Coloring and Proportional Coloring}

The study of equitable coloring began with a conjecture of Erd\H{o}s in 1964~\cite{E64}, although the concept was not formally introduced until Meyer's work in 1973~\cite{M73}. An \emph{equitable $k$-coloring} of a graph $G$ is a proper $k$-coloring $f$ such that the sizes of any two color classes differ by at most one. Here, we adopt the convention that a $k$-coloring has exactly $k$ color classes, some of which may be empty. Consequently, each color class in an equitable $k$-coloring has size either $\lceil |V(G)|/k \rceil$ or $\lfloor |V(G)|/k \rfloor$. We say that a graph $G$ is \emph{equitably $k$-colorable} if it admits an equitable $k$-coloring.

List coloring was introduced independently by Vizing~\cite{V76} and Erd\H{o}s, Rubin, and Taylor~\cite{ER80} in the 1970s. A \emph{list assignment} $L$ for a graph $G$ is a function that assigns to each vertex of $G$ a set of colors. If $|L(v)|=k$ for every $v \in V(G)$, then $L$ is called a \emph{$k$-assignment} for $G$. A graph $G$ is \emph{$L$-colorable} if there exists a proper coloring $f$ of $G$ such that $f(v) \in L(v)$ for every $v \in V(G)$. Such a coloring is called a \emph{proper $L$-coloring} of $G$. We say that $G$ is \emph{$k$-choosable} if it admits a proper $L$-coloring for every $k$-assignment $L$ of $G$. Both list coloring and equitable coloring have been extensively studied since their introduction.

In 2003, Kostochka, Pelsmajer, and West introduced a list analogue of equitable coloring, called \emph{equitable choosability}~\cite{KP03}. Suppose that $L$ is a $k$-assignment for a graph $G$. A proper $L$-coloring of $G$ is \emph{equitable} if each color appears on at most $\lceil |V(G)|/k \rceil$ vertices. Such a coloring is called an \emph{equitable $L$-coloring} of $G$, and we say that $G$ is \emph{equitably $L$-colorable} if it admits an equitable $L$-coloring. We say that $G$ is \emph{equitably $k$-choosable} if $G$ is equitably $L$-colorable for every $k$-assignment $L$ of $G$. Observe that equitable choosability imposes the same upper bound on color class sizes as equitable coloring, but it does not require the sizes of the color classes to be balanced.

In 2019, another list analogue of equitable coloring, called \emph{proportional choosability}~\cite{KM19}, was introduced. It imposes both upper and lower bounds on color class sizes based on the availability of each color. Since a color may appear in as few as one list, we cannot require every color class to have size at least $\lfloor |V(G)|/k \rfloor$. Instead, proportional choosability uses bounds that depend on the number of lists containing each color.

Let $L$ be a list assignment for a graph $G$. The \emph{palette of colors associated with $L$} is the set $\bigcup_{v \in V(G)} L(v)$, which we denote by $\mathcal{L}$ throughout this paper. For each color $c \in \mathcal{L}$, the \emph{multiplicity of $c$ in $L$}, denoted by $\mu_L(c)$ (or simply $\mu(c)$ when $L$ is clear from context), is the number of lists in which $c$ appears. In particular, $\mu_L(c) = |\{v \in V(G) : c \in L(v)\}|$.

Suppose that $L$ is a $k$-assignment for a graph $G$. A proper $L$-coloring $f$ of $G$ is called a \emph{proportional $L$-coloring} if, for every $c \in \mathcal{L}$, the color class $f^{-1}(c)$ has size either $\lfloor \mu_L(c)/k \rfloor$ or $\lceil \mu_L(c)/k \rceil$. We say that $G$ is \emph{proportionally $L$-colorable} if it admits a proportional $L$-coloring. Furthermore, we say that $G$ is \emph{proportionally $k$-choosable} if it is proportionally $L$-colorable for every $k$-assignment $L$ of $G$.

We now recall several basic properties of proportional choosability. For the first property, recall that if a graph is $k$-choosable, then it is $k$-colorable. However, the analogous implication fails for equitable choosability: a graph may be equitably $k$-choosable without being equitably $k$-colorable.

\begin{pro}[\cite{KM19}] \label{pro: motivation}
If $G$ is proportionally $k$-choosable, then $G$ is both equitably $k$-choosable and equitably $k$-colorable.
\end{pro}

The next two properties highlight additional advantages of proportional choosability. In particular, proportional choosability is monotone in $k$ and is also preserved under taking subgraphs.

\begin{pro}[\cite{KM19}] \label{pro: monoink}
If $G$ is proportionally $k$-choosable, then $G$ is proportionally $(k+1)$-choosable.
\end{pro}

\begin{pro}[\cite{KM19}] \label{lem: monotone}
Suppose $H$ is a subgraph of $G$. If $G$ is proportionally $k$-choosable, then $H$ is proportionally $k$-choosable.
\end{pro}

For any graph $G$, the \emph{proportional choice number} of $G$, denoted by $\chi_{pc}(G)$, is the smallest positive integer $k$ such that $G$ is proportionally $k$-choosable. By Proposition~\ref{pro: monoink}, a graph $G$ is proportionally $k$-choosable if and only if $k \geq \chi_{pc}(G)$.

We conclude by recalling a characterization of proportionally $2$-choosable graphs. Recall that a \emph{linear forest} is a graph whose components are paths.

\begin{thm}[\cite{KM20}] \label{thm: prop2}
A graph $G$ is proportionally $2$-choosable if and only if $G$ is a linear forest whose largest component has at most five vertices and all other components have at most two vertices.
\end{thm}

\subsection{Paper Focus}

In~\cite{M26}, the second named author showed that $\chi_{pc}(P_n)\leq 4$ for every $n\in\N$, answering the question of whether the proportional choice numbers of paths are bounded by a constant. The same paper also showed that the inequality $\chi_{pc}(G)\leq\Delta(G)+1$ does not hold for all graphs. These results motivated the following two conjectures which are the focus of this paper.

\begin{conj} \label{conj: its3}
For every $n \geq 6$, $\chi_{pc}(P_n) = 3$.
\end{conj}

\begin{conj} \label{conj: upper}
For every graph $G$ with at least one edge, $\chi_{pc}(G) \leq \Delta(G) + \col(G) - 1$.
\end{conj}

Notice that Conjecture~\ref{conj: upper} would immediately imply Conjecture~\ref{conj: its3}. Indeed, when $n\geq 6$, we know that $\Delta(P_n)=2$ and $\col(P_n)=2$, so Conjecture~\ref{conj: upper} would imply $\chi_{pc}(P_n)\leq 3$. Moreover, Theorem~\ref{thm: prop2} implies that $\chi_{pc}(P_n)\geq 3$ whenever $n\geq 6$. Thus, Conjecture~\ref{conj: upper} would imply that $\chi_{pc}(P_n)=3$ for every $n\geq 6$.

It is also worth mentioning that the upper bound proposed in Conjecture~\ref{conj: upper} is best possible. Indeed, it was shown in~\cite{M26} that if $a,b\in\N$ and $G=(a+b-1)K_{a,b}$, then $\chi_{pc}(G)=a+b$. Since $\Delta(G)=\max\{a,b\}$ and $\col(G)=\min\{a,b\}+1$, we have $\chi_{pc}(G)=\Delta(G)+\col(G)-1$. Consequently, equality is attained for infinitely many graphs.

The remainder of this paper is organized as follows. In Section~\ref{generalreduction}, we introduce the necessary background on abstract simplicial complexes and matroids and establish a lower bound on the connectivity of the independence complex of a graph. In Section~\ref{main}, we use these results together with a theorem of Aharoni and Berger to prove Conjecture~\ref{conj: upper} which consequently resolves Conjecture~\ref{conj: its3} as well.

\section{Simplicial Complex and Matroid Preliminaries}\label{generalreduction}

\subsection{Abstract Simplicial Complexes}

Recall that an \emph{abstract simplicial complex}, $\mathcal{C}$, on the finite ground set $V$ is a nonempty set of subsets of $V$ with the property that if $A \in \mathcal{C}$, then every subset of $A$ is an element of $\mathcal{C}$. The elements of $\mathcal{C}$ are called the \emph{simplices} of $\mathcal{C}$.

Suppose $\mathcal{C}$ is an abstract simplicial complex on the ground set $V$ and $X \subseteq V$. We let $\mathcal{C}|X$ denote the abstract simplicial complex consisting of all simplices of $\mathcal{C}$ that are subsets of $X$. That is,
\[
\mathcal{C}|X = \{A \in \mathcal{C} : A \subseteq X\}.
\]
We call $\mathcal{C}|X$ the \emph{induced subcomplex} of $\mathcal{C}$ on $X$.

It is well-known that every abstract simplicial complex on a finite ground set has a geometric realization that is unique up to homeomorphism. Indeed, suppose $V = \{v_1, \dots, v_n\}$ is the ground set of the abstract simplicial complex $\mathcal{C}$. Then, let $e_1, \dots, e_n$ be the standard basis vectors of $\mathbb{R}^n$. For every nonempty simplex $A \in \mathcal{C}$, its geometric realization, denoted by $|A|$, is the convex hull of $\{e_i : v_i \in A\}$. Then the geometric realization of $\mathcal{C}$ is
\[
\|\mathcal{C}\| = \bigcup_{\substack{S \in \mathcal{C}\\S \neq \emptyset}} |S|,
\]
equipped with the subspace topology induced by the Euclidean metric on $\mathbb{R}^n$.

We say that an abstract simplicial complex $\mathcal{C}$ is \emph{$k$-connected} for a nonnegative integer $k$ if $\|\mathcal{C}\|$ is nonempty and, for every integer $0 \leq i \leq k$, every continuous map from the $i$-dimensional sphere $S^i$ into $\|\mathcal{C}\|$ can be extended to a continuous map from the $(i+1)$-dimensional ball $B^{i+1}$ into $\|\mathcal{C}\|$. By convention, $\mathcal{C}$ is $(-1)$-connected if $\|\mathcal{C}\|$ is nonempty. Following Aharoni and Berger~\cite{AB06}, we define $\eta(\mathcal{C})$ to be the largest integer $k$ for which $\mathcal{C}$ is $k$-connected, plus two. If $\mathcal{C}$ is $k$-connected for every $k$, then we set $\eta(\mathcal{C})=\infty$. Finally, if $\|\mathcal{C}\|$ is empty, we set $\eta(\mathcal{C})=0$.

Suppose $\mathcal{C}_1,\ldots,\mathcal{C}_m$ are abstract simplicial complexes on pairwise disjoint ground sets $V_1,\ldots,V_m$. We define the \emph{join} of these complexes to be the abstract simplicial complex
\[
\bigvee_{i=1}^m \mathcal{C}_i = \left\{A_1\cup\cdots\cup A_m : A_i\in\mathcal{C}_i \text{ for each } i\in[m]\right\}.
\]
The following observation follows by repeatedly applying the result of Aharoni and Berger~\cite[Lemma~2.2]{AB06}.

\begin{obs} \label{obs: join}
Suppose $\mathcal{C}_1,\ldots,\mathcal{C}_m$ are abstract simplicial complexes on pairwise disjoint ground sets. Then,
\[
\eta\left(\bigvee_{i=1}^m\mathcal{C}_i\right) \geq \sum_{i=1}^m\eta(\mathcal{C}_i).
\]
\end{obs}

Suppose $\mathcal{C}$ is an abstract simplicial complex with ground set $V$. Suppose there is a $v \in V$ such that, for each $A \in \mathcal{C}$, we have $(A \cup \{v\}) \in \mathcal{C}$. When such a $v$ exists, we say that $\mathcal{C}$ is a \emph{cone} and $v$ is an \emph{apex vertex} of $\mathcal{C}$. The following observation is easy to prove.

\begin{obs} \label{obs: cone}
Suppose the abstract simplicial complex $\mathcal{C}$ is a cone with ground set $V$ and apex vertex $v$. Then, $\eta(\mathcal{C})=\infty$.
\end{obs}

\begin{proof}
Let $k$ be an arbitrary nonnegative integer. Suppose $f: S^k \rightarrow \|\mathcal{C}\|$ is an arbitrary continuous map from the $k$-dimensional sphere into the geometric realization of $\mathcal{C}$. In addition, suppose that $v$ is associated with $e_1$ in the geometric realization of $\mathcal{C}$.

Now, let $g: B^{k+1} \rightarrow \|\mathcal{C}\|$ be the function given by
\[
g(\boldsymbol{x}) = \begin{cases}
    (1 - \|\boldsymbol{x}\|)e_1 + \|\boldsymbol{x}\| f \left( \frac{\boldsymbol{x}}{\|\boldsymbol{x}\|} \right) & \text{if } \boldsymbol{x} \neq \boldsymbol{0}, \\
    e_1 & \text{if } \boldsymbol{x} = \boldsymbol{0}.
\end{cases}
\]
Since $\mathcal{C}$ is a cone with apex vertex $v$, the line segment joining $e_1$ to any point in $\|\mathcal{C}\|$ is contained in $\|\mathcal{C}\|$. Thus, the outputs of $g$ are indeed in $\|\mathcal{C}\|$. It is also easy to check that $g$ is continuous, including at $\boldsymbol{0}$. Moreover, $g(\boldsymbol{x})=f(\boldsymbol{x})$ whenever $\boldsymbol{x}\in S^k$ which means $g$ is an extension of $f$. Since $k$ and $f$ were arbitrary, the observation follows.
\end{proof}

We now turn our attention to the abstract simplicial complexes that play a crucial role in our proof of Conjecture~\ref{conj: upper}. Suppose that $G$ is a graph. We let $\mathcal{I}(G)$ be the set containing all independent sets of vertices of $G$. It is easy to check that $\mathcal{I}(G)$ is an abstract simplicial complex with ground set $V(G)$. Furthermore, if $v \in V(G)$ is an isolated vertex in $G$, then $\mathcal{I}(G)$ is a cone with apex vertex $v$. Consequently, Observation~\ref{obs: cone} implies that $\eta(\mathcal{I}(G)) = \infty$ whenever $G$ has an isolated vertex.  Also, note that if $X \subseteq V(G)$, then $\mathcal{I}(G)|X = \mathcal{I}(G[X])$, where $G[X]$ is the subgraph of $G$ induced by $X$.

The following observation follows from a standard result on the connectivity of CW complexes and their subcomplexes (see, for example, \cite[Corollary~4.12]{Hat02}).

\begin{obs} \label{obs: skeleton}
Suppose $t \in \N$ and $\mathcal{C}$ is an abstract simplicial complex whose largest simplices have $t+1$ elements, and let $a$ be an integer with $2 \leq a \leq t+1$. Let $\mathcal{C}'$ be the abstract simplicial complex obtained from $\mathcal{C}$ by removing all simplices having between $a$ and $t+1$ elements. Then,
\[
\eta(\mathcal{C}') \geq \min\{\eta(\mathcal{C}),a-1\}.
\]
\end{obs}

The following theorem of Meshulam~\cite{Me03}, stated in terms of homotopy connectivity by Haxell and Wdowinski~\cite[Theorem~13]{HW24}, provides a useful recursive lower bound on $\eta(\mathcal{I}(G))$ that we will use in our proof of Conjecture~\ref{conj: upper}. Before stating the theorem, we introduce a graph operation. Suppose $e=uv$ is an edge in a graph $G$. We let $G*e$ denote the graph obtained from $G$ by deleting all vertices in $N_G[u]\cup N_G[v]$. The operation of obtaining $G*e$ from $G$ is called \emph{exploding} the edge $e$.

\begin{thm}[Meshulam's Theorem] \label{thm: meshulam}
Suppose $G$ is a graph and $e\in E(G)$. Then,
\[
\eta(\mathcal{I}(G)) \geq \min\left\{\eta(\mathcal{I}(G-e)), \eta(\mathcal{I}(G*e))+1\right\}.
\]
\end{thm}

We are now ready to prove a key lemma.

\begin{lem} \label{lem: independent}
Suppose $G$ is a graph on $n$ vertices with at least one edge. Then,
\[
\eta(\mathcal{I}(G)) \geq \frac{n}{\Delta(G)+\col(G)-1}.
\]
\end{lem}

\begin{proof}
The proof is by induction on $k = n + |E(G)|$. In the case $G = K_2$, the desired result follows since $\|\mathcal{I}(G)\|$ is nonempty and $\Delta(G)+\col(G)-1=2$.

So, suppose that $k \geq 4$, $|E(G)| \geq 1$, and the desired result holds for all graphs with at least one edge for which the sum of the number of vertices and edges is less than $k$. We may suppose that $G$ has no isolated vertices since we know that when $G$ has an isolated vertex $\mathcal{I}(G)$ is a cone and the desired result holds by Observation~\ref{obs: cone}.

Now, suppose that $x$ is a vertex in $G$ of minimum degree and suppose that $e=xy \in E(G)$. Note that $d_G(x) \leq \col(G)-1$ and $d_G(y) \leq \Delta(G)$. Let $G_1 = G-e$ and $G_2 = G*e$. Notice that $G_1$ has coloring number at most $\col(G)$ and maximum degree at most $\Delta(G)$. The same holds for $G_2$ whenever $G_2$ is nonempty, since deleting edges or vertices of a graph cannot increase these parameters.

Let $D=\Delta(G)+\col(G)-1$. First, suppose that $G_2$ has no vertices. Since $G_2$ is obtained by exploding the edge $xy$, we know that $n \leq d_G(x)+d_G(y)\leq D$. Since $\|\mathcal{I}(G)\|$ is nonempty, $\eta(\mathcal{I}(G))\geq 1\geq n/D$, and the desired result follows. So, we may suppose that $G_2$ has at least one vertex.

If $G_1$ has an isolated vertex, then $\eta(\mathcal{I}(G_1))=\infty$ by Observation~\ref{obs: cone}. Otherwise, the induction hypothesis implies
\[
\eta(\mathcal{I}(G_1)) \geq \frac{|V(G_1)|}{\Delta(G_1)+\col(G_1)-1} \geq \frac{n}{D}.
\]

Similarly, if $G_2$ has an isolated vertex, then $\eta(\mathcal{I}(G_2))=\infty$ by Observation~\ref{obs: cone}. Otherwise, note that $G_2$ has at least $n-d_G(x)-d_G(y)$ vertices. So, by the induction hypothesis and the facts we have observed,
\begin{align*}
\eta(\mathcal{I}(G_2))+1 &\geq \frac{|V(G_2)|}{\Delta(G_2)+\col(G_2)-1}+1 \\
&\geq \frac{n-d_G(x)-d_G(y)}{D}+1 \\
&= \frac{n+(\col(G)-1-d_G(x))+(\Delta(G)-d_G(y))}{D} \\
&\geq \frac{n}{D}.
\end{align*}
Thus, in all cases, $\eta(\mathcal{I}(G_1))\geq n/D$ and $\eta(\mathcal{I}(G_2))+1\geq n/D$. Meshulam's Theorem then implies the desired inequality for $\mathcal{I}(G)$ and the induction step is complete.
\end{proof}

\subsection{Matroids}

We now give a brief introduction to matroid theory, which will give the necessary background for Theorem~\ref{thm: ABftw}, as well as the proof of Lemma \ref{lem: key}. We do not assume any prior knowledge of matroids.

A \emph{matroid} is an ordered pair $\mathcal{M}=(V,\mathcal{I})$, where $V$ is a finite set, called the \emph{ground set}, and $\mathcal{I}$ is a collection of subsets of $V$, called \emph{independent sets}, such that \begin{enumerate}
\item $\mathcal{I}$ is a nonempty abstract simplicial complex and
\item if $I_1$ and $I_2$ are members of $\mathcal{I}$ with $|I_2|>|I_1|$, then there is an element $v\in I_2-I_1$ such that $I_1\cup\{v\}\in\mathcal{I}$.
\end{enumerate}

Since our main result relies on Theorem~\ref{thm: ABftw}, proved by Aharoni and Berger in \cite{AB06}, we largely follow the notation of \cite{AB06}. However, a standard reference work in matroid theory is Oxley \cite{o11}. It is perhaps worth pointing out some notational differences between \cite{AB06} and \cite{o11}.

In \cite{o11} and many other works on matroid theory, the ground set of a matroid is denoted by $E$. This is in part because a fundamental example of a matroid is the cycle matroid of a graph, whose ground set is the edge set of the graph. Much standard matroid terminology comes from this generalization of graphs. In this paper, however, we generalize graphs in a different way. A simple graph can be thought of as an abstract simplicial complex whose ground set is the vertex set of the graph and whose nonempty simplices have size 2 (edges) or size 1 (vertices).

An implication of the definition of a matroid is that the inclusion-wise maximal independent sets of a matroid $\mathcal{M}$ all have the same size. These maximal independent sets are the \emph{bases} of the matroid, and their size is the \emph{rank} of the matroid, denoted by $\rho(\mathcal{M})$. Additionally, every matroid $\mathcal{M}$ on ground set $V$ has a rank function $\rho_{\mathcal{M}}$ such that $\rho_{\mathcal{M}}(X)$ is the size of the largest independent subset of $X$, for each subset $X\subseteq V$.

It follows from the definition of the rank function given above that if $X\subseteq Y\subseteq V$, then $\rho_{\mathcal{M}}(X)\leq\rho_{\mathcal{M}}(Y)$. The rank-maximal sets in a matroid are called the $\emph{flats}$ of the matroid. In other words, a flat of a matroid $\mathcal{M}$ on ground set $V$ is a set $F\subseteq V$ such that $\rho_{\mathcal{M}}(F\cup\{v\})>\rho_{\mathcal{M}}(F)$ for each element $v\in V-F$.

The specific matroids that we need in this paper have a relatively simple structure. Let $V=V_1\cup V_2\cup\cdots V_n$ be a disjoint union of nonempty sets. The \emph{partition matroid} associated with this partition has as its collection of independent sets $\mathcal{I}=\{X\subseteq V: |X\cap V_i|\leq1\text{ for each }i\in[n]\}$. It is not difficult to check that a partition matroid is indeed a matroid. Note that the rank function of a partition matroid assigns to each subset $X\subseteq V$ the number of $V_i$ with which $X$ has nonempty intersection. Moreover, the flats of a partition matroid are the sets of the form $\bigcup_{i\in S}V_i$, where $S\subseteq [n]$.

\section{Proving Conjecture~\ref{conj: upper}} \label{main}

The following theorem of Aharoni and Berger will be key to our proof.

\begin{thm}[\cite{AB06}, Theorem~4.5] \label{thm: ABftw}
Suppose $\mathcal{M}$ is a matroid and $\mathcal{C}$ is an abstract simplicial complex on the same ground set $V$. If for every flat $F$ of $\mathcal{M}$,
\[
\eta(\mathcal{C}|(V - F)) \geq \rho(\mathcal{M})-\rho_{\mathcal{M}}(F),
\]
then there is a basis of $\mathcal{M}$ in $\mathcal{C}$.
\end{thm}

We now introduce some notions and notation we will need. Suppose that $G$ is a graph and $L$ is a $k$-assignment for $G$ such that for each $c \in \mathcal{L}$, $k$ divides $\mu_L(c)$. We let $\mathcal{P}(G,L) = \{(v,c) : v \in V(G), c \in L(v)\}$, and for each $v \in V(G)$, we let $P_v(G,L) = \{(v,c) : c \in L(v)\}$.

We use $\mathcal{M}(G,L)$ to denote the partition matroid on the ground set $\mathcal{P}(G,L)$ with partition $\{P_v(G,L) : v \in V(G)\}$, where a set is independent if it contains at most one element from each $P_v(G,L)$. Note that each basis of $\mathcal{M}(G,L)$ corresponds to an $L$-coloring (not necessarily proper) of $G$.

For each $c \in \mathcal{L}$, we let $\mathcal{C}_c(G,L)$ denote the abstract simplicial complex with ground set $Q_c(G,L)=\{(v,c) : v \in V(G), c \in L(v)\}$ whose simplices are all sets of the form $\{(v,c) : v \in I\}$, where $I \in \mathcal{I}(G)$, $c \in L(v)$ for every $v \in I$, and $|I| \leq \mu_L(c)/k$. Notice that $|Q_c(G,L)|=\mu_L(c)$. Furthermore, if $G_c = G[\{v \in V(G) : c \in L(v)\}]$, then $\mathcal{C}_c(G,L)$ is isomorphic to the abstract simplicial complex formed from $\mathcal{I}(G_c)$ by removing all simplices with more than $\mu_L(c)/k$ elements.

Since $\{Q_c(G,L) : c \in \mathcal{L}\}$ is a partition of $\mathcal{P}(G,L)$, the abstract simplicial complex $\mathcal{C}(G,L)$ given by $\mathcal{C}(G,L) = \bigvee_{c \in \mathcal{L}} \mathcal{C}_c(G,L)$ is an abstract simplicial complex with ground set $\mathcal{P}(G,L)$.

For all the notions introduced above, we may drop $(G,L)$ when the graph and list assignment are clear from context. We are now ready for our main lemma that will allow us to prove Conjecture~\ref{conj: upper}.

\begin{lem} \label{lem: key}
Suppose $G$ is a graph with at least one edge and $k = \Delta(G)+\col(G)-1$. If $L$ is a $k$-assignment for $G$ such that for each $c \in \mathcal{L}$, $k$ divides $\mu_L(c)$, then there is a proportional $L$-coloring of $G$.
\end{lem}

\begin{proof}
Let $n = |V(G)|$. Consider the matroid $\mathcal{M}(G,L)$ and the abstract simplicial complex $\mathcal{C}(G,L)$ on the ground set $\mathcal{P}(G,L)$.

Suppose that $F$ is an arbitrary flat of $\mathcal{M}$. Then, there is some $S \subseteq V(G)$ such that $F = \bigcup_{v \in S} P_v$. We wish to apply Theorem~\ref{thm: ABftw}. So, we must show
\[
\eta(\mathcal{C}|(\mathcal{P} - F)) \geq \rho(\mathcal{M})-\rho_{\mathcal{M}}(F).
\]
We first note that $\rho(\mathcal{M})-\rho_{\mathcal{M}}(F) = n - |S|$. In the case $S = V(G)$, $\mathcal{P}-F = \emptyset$, and the desired inequality clearly holds. So, we may assume that $S \subsetneq V(G)$.

For each $c \in \mathcal{L}$, let $\mathcal{C}_c' = \mathcal{C}_c|(Q_c - F)$ and let $G_c' = G[\{v \in V(G)-S : c \in L(v)\}]$. In the case $\mathcal{C}_c'$ contains only the empty set, we have $\eta(\mathcal{C}_c')=0=|V(G_c')|/k$. So, suppose that $G_c'$ has at least one vertex. Notice that $\mathcal{C}_c'$ is isomorphic to the abstract simplicial complex formed from $\mathcal{I}(G_c')$ by removing all simplices with more than $\mu_L(c)/k$ elements. So, by Observation~\ref{obs: skeleton}, or trivially if no simplices are removed, we have
\[
\eta(\mathcal{C}_c') \geq \min\left\{\eta(\mathcal{I}(G_c')), \frac{\mu_L(c)}{k}\right\}.
\]
In the case $G_c'$ has no edges, $\mathcal{I}(G_c')$ is a cone, and Observation~\ref{obs: cone} implies $\eta(\mathcal{I}(G_c'))=\infty$. Otherwise, since $\Delta(G_c')+\col(G_c')-1 \leq k$, Lemma~\ref{lem: independent} implies
\[
\eta(\mathcal{I}(G_c')) \geq \frac{|V(G_c')|}{\Delta(G_c')+\col(G_c')-1} \geq \frac{|V(G_c')|}{k}.
\]
Since $\mu_L(c)\geq |V(G_c')|$, we therefore have
\[
\eta(\mathcal{C}_c') \geq \min\left\{\frac{|V(G_c')|}{k},\frac{\mu_L(c)}{k}\right\} = \frac{|V(G_c')|}{k}.
\]
Thus, in all cases, even when $G_c'$ has no vertices, we have $\eta(\mathcal{C}_c') \geq |V(G_c')|/k$.

Now, let $L'$ be the list assignment for $G-S$ obtained by restricting the domain of $L$ to $V(G)-S$.  We extend the multiplicity function $\mu_{L'}$ to all of $\mathcal{L}$ by setting $\mu_{L'}(c)=0$ whenever $c$ is not in the palette of $L'$. Since $|V(G_c')| = \mu_{L'}(c)$ and $\mathcal{C}|(\mathcal{P}-F) = \bigvee_{c \in \mathcal{L}} \mathcal{C}_c'$, Observation~\ref{obs: join} implies
\[
\eta(\mathcal{C}|(\mathcal{P} - F)) \geq \sum_{c \in \mathcal{L}} \eta(\mathcal{C}_c') \geq \sum_{c \in \mathcal{L}} \frac{|V(G_c')|}{k} = \frac{1}{k} \sum_{c \in \mathcal{L}} \mu_{L'}(c) = \frac{1}{k} k(n-|S|) = n-|S|
\]
as desired. Consequently, Theorem~\ref{thm: ABftw} implies there is a basis $B$ of $\mathcal{M}$ in $\mathcal{C}$.

Since $B$ is a basis of $\mathcal{M}$, it corresponds to an $L$-coloring $f$ of $G$. Furthermore, since $B \in \mathcal{C}$, the coloring $f$ is proper and $|f^{-1}(c)|\leq \mu_L(c)/k$ for every $c \in \mathcal{L}$. Since
\[
\sum_{c \in \mathcal{L}} \frac{\mu_L(c)}{k} = n = \sum_{c \in \mathcal{L}} |f^{-1}(c)|,
\]
we must have $|f^{-1}(c)|=\mu_L(c)/k$ for every $c \in \mathcal{L}$. Thus, $f$ is a proportional $L$-coloring of $G$, and the proof is complete.
\end{proof}

Finally, we prove Conjecture~\ref{conj: upper}.

\begin{thm} \label{thm: 9years}
For any graph $G$ with at least one edge, $\chi_{pc}(G) \leq \Delta(G) + \col(G)-1$.
\end{thm}

\begin{proof}
Suppose $G$ has $n$ vertices, and let $k = \Delta(G)+\col(G)-1$. Since $G$ has an edge, $k \geq 2$. Suppose that $L$ is an arbitrary $k$-assignment for $G$. To prove the desired result, we must show that there is a proportional $L$-coloring of $G$. By relabeling colors if necessary, we may assume that $\N^2$ is disjoint from $\mathcal{L}$, and we let $z$ be a color not in $\mathcal{L}\cup\N^2$.

Let $C = \{ c \in \mathcal{L} : \text{$\mu_L(c)$ is not divisible by $k$}\}$. If $C = \emptyset$, a proportional $L$-coloring of $G$ exists by Lemma~\ref{lem: key}. So, we may assume that $C \neq \emptyset$ and the elements of $C$ are named so that $C = \{c_1, \dots, c_t\}$. For each $i \in [t]$, associate with $c_i$ a set of $k$ vertices $C_i = \{(c_i,1), \dots, (c_i,k)\}$ and assume that $C_i$ is disjoint from $V(G)$.

Now, let $G'$ be the graph with $V(G') = V(G) \cup \bigcup_{i=1}^t C_i$ and $E(G')=E(G)$. Notice that $\Delta(G') = \Delta(G)$ and $\col(G') = \col(G)$. Now, we will construct a $k$-assignment $L'$ for $G'$ such that we can apply Lemma~\ref{lem: key} to $G'$ and $L'$.

For each $i \in [t]$, let $r_i = (\mu_L(c_i) \mod k)$ and $d_i = k-r_i$. Note that since $\sum_{c \in \mathcal{L}} \mu_L(c) = nk$, there is a $q \in \N$ such that $\sum_{i=1}^t r_i = qk$. Now, let $L'(v) = L(v)$ for each $v \in V(G)$. Also, for each $i \in [t]$, let $L'((c_i,j)) = \{c_i\} \cup \{(i,p) : p \in [k-1]\}$ if $j \in [d_i]$ and let $L'((c_i,j)) = \{z\} \cup \{(i,p) : p \in [k-1]\}$ if $j \in \{d_i+1, \dots, k\}$.

Notice that for each $c \in C$, $\mu_{L'}(c) = k \lceil \mu_L(c)/k \rceil$, and for each $c \in \mathcal{L}-C$, $\mu_{L'}(c)=\mu_L(c)$. Also, for each $(i,j) \in [t] \times [k-1]$, $\mu_{L'}((i,j)) = k$. Finally, $\mu_{L'}(z) = qk$. So, by Lemma~\ref{lem: key}, there is a proportional $L'$-coloring, $f'$, of $G'$. If $f$ is the restriction of $f'$ to $V(G)$, we claim that $f$ is a proportional $L$-coloring of $G$.

To see this, first notice that $f$ is clearly proper. Next, any $c \in \mathcal{L}-C$ is used exactly $\mu_L(c)/k$ times by $f'$ exclusively on vertices in $V(G)$ since such a color $c$ only appears in lists associated with vertices in $V(G)$ by $L'$. Then, for each $i \in [t]$, notice that $f'$ uses $c_i$ exactly $\lceil \mu_L(c_i)/k \rceil$ times. Also, $c_i$ only appears in lists associated with vertices in $V(G) \cup \{(c_i,j) : j \in [d_i]\}$ by $L'$. Since $\mu_{L'}((i,j)) = k$ for each $j \in [k-1]$, each color in $\{(i,j): j \in [k-1]\}$ must be used exactly once by $f'$ on a vertex in $C_i$. Thus, exactly $k-1$ vertices in $C_i$ receive colors from this set, and the remaining vertex receives either $c_i$ or $z$. Consequently, $c_i$ is used at most once by $f'$ on a vertex in $C_i$, and it immediately follows that $f$ uses $c_i$ either $\lceil \mu_L(c_i)/k \rceil$ or $\lfloor \mu_L(c_i)/k \rfloor$ times.

So, we have that $f$ is a proportional $L$-coloring of $G$, and our proof is complete.
\end{proof}

We conclude by determining the proportional choice number of any path, thereby resolving Conjecture~\ref{conj: its3}.

\begin{cor} \label{cor: paths}
Suppose $n \in \N$. Then, $\chi_{pc}(P_1)=1$, $\chi_{pc}(P_n)=2$ when $2\leq n\leq 5$, and $\chi_{pc}(P_n)=3$ when $n\geq 6$.
\end{cor}

\begin{proof}
Clearly, $\chi_{pc}(P_1)=1$. Since $P_n$ has an edge whenever $n\geq 2$, we know that $\chi_{pc}(P_n)\geq 2$. By Theorem~\ref{thm: prop2}, $\chi_{pc}(P_n)=2$ when $2\leq n\leq 5$ and $\chi_{pc}(P_n)\geq 3$ when $n\geq 6$. Finally, when $n\geq 6$, $\Delta(P_n)=2$ and $\col(P_n)=2$. So, Theorem~\ref{thm: 9years} implies $\chi_{pc}(P_n)\leq 3$, and the result follows.
\end{proof}

\par \medskip \textbf{Use of generative AI.} During the preparation of this manuscript, the authors used ChatGPT (OpenAI) and Claude (Anthropic) as supplementary tools to examine technical details, identify applicable mathematical techniques, and assist with editing and presentation. All mathematical arguments were independently reviewed and verified by the authors, who take full responsibility for the content of the manuscript.

\bibliographystyle{hplain}
\bibliography{bibliography}{}

@article{AB06,
  author  = {Aharoni, R. and Berger, E.},
  title   = {The intersection of a matroid and a simplicial complex},
  journal = {Trans. Amer. Math. Soc.},
  volume  = {358},
  number  = {11},
  year    = {2006},
  pages   = {4895--4917},
  doi     = {10.1090/S0002-9947-06-03833-5}
}

@incollection{E64,
  author    = {Erd{\H{o}}s, P.},
  title     = {Problem 9},
  booktitle = {Theory of Graphs and Its Applications},
  editor    = {Fiedler, M.},
  publisher = {Publ. House Czechoslovak Acad. Sci.},
  address   = {Prague},
  year      = {1964},
  pages     = {159},
  note      = {Proc. Sympos., Smolenice, 1963}
}

@inproceedings{ER80,
  author    = {Erd{\H{o}}s, P. and Rubin, A. L. and Taylor, H.},
  title     = {Choosability in graphs},
  booktitle = {Proceedings of the West Coast Conference on Combinatorics, Graph Theory and Computing (Humboldt State Univ., Arcata, Calif., 1979)},
  series    = {Congr. Numer.},
  volume    = {XXVI},
  publisher = {Utilitas Math.},
  address   = {Winnipeg, MB},
  year      = {1980},
  pages     = {125--157}
}

@book{Hat02,
  author    = {Hatcher, A.},
  title     = {Algebraic Topology},
  publisher = {Cambridge University Press},
  address   = {Cambridge},
  year      = {2002}
}

@article{HW24,
  author  = {Haxell, P. and Wdowinski, R.},
  title   = {Degree criteria and stability for independent transversals},
  journal = {J. Graph Theory},
  volume  = {106},
  number  = {2},
  year    = {2024},
  pages   = {352--371},
  doi     = {10.1002/jgt.23085}
}

@article{KM19,
  author  = {Kaul, H. and Mudrock, J. A. and Pelsmajer, M. J. and Reiniger, B.},
  title   = {Proportional choosability: a new list analogue of equitable coloring},
  journal = {Discrete Math.},
  volume  = {342},
  number  = {8},
  year    = {2019},
  pages   = {2371--2383},
  doi     = {10.1016/j.disc.2019.05.011}
}

@article{KM20,
  author  = {Kaul, H. and Mudrock, J. A. and Pelsmajer, M. J. and Reiniger, B.},
  title   = {A simple characterization of proportionally {$2$}-choosable graphs},
  journal = {Graphs Combin.},
  volume  = {36},
  number  = {3},
  year    = {2020},
  pages   = {679--687},
  doi     = {10.1007/s00373-020-02146-z}
}

@article{KP03,
  author  = {Kostochka, A. V. and Pelsmajer, M. J. and West, D. B.},
  title   = {A list analogue of equitable coloring},
  journal = {J. Graph Theory},
  volume  = {44},
  number  = {3},
  year    = {2003},
  pages   = {166--177},
  doi     = {10.1002/jgt.10137}
}

@article{Me03,
  author  = {Meshulam, R.},
  title   = {Domination numbers and homology},
  journal = {J. Combin. Theory Ser. A},
  volume  = {102},
  number  = {2},
  year    = {2003},
  pages   = {321--330},
  doi     = {10.1016/S0097-3165(03)00045-1}
}

@article{M73,
  author  = {Meyer, W.},
  title   = {Equitable coloring},
  journal = {Amer. Math. Monthly},
  volume  = {80},
  number  = {8},
  year    = {1973},
  pages   = {920--922},
  doi     = {10.2307/2319405}
}

@misc{M26,
  author        = {Mudrock, J. A.},
  title         = {Improved bounds on proportional choice numbers},
  year          = {2026},
  eprint        = {2610.33914},
  archivePrefix = {arXiv},
  primaryClass  = {math.CO},
  doi           = {10.48550/arXiv.2610.33914},
  url           = {https://arxiv.org/abs/2610.33914}
}

@article{V76,
  author  = {Vizing, V. G.},
  title   = {Coloring the vertices of a graph in prescribed colors},
  journal = {Diskret. Analiz},
  number  = {29},
  year    = {1976},
  pages   = {3--10},
  note    = {Metody Diskret. Anal. v Teorii Kodov i Shem, 101 (in Russian)}
}

@book{W01,
  author    = {West, D. B.},
  title     = {Introduction to Graph Theory},
  edition   = {2nd},
  publisher = {Prentice Hall},
  address   = {Upper Saddle River, NJ},
  year      = {2001}
}

@book {o11,
    AUTHOR = {Oxley, James},
     TITLE = {Matroid Theory},
    SERIES = {Oxford Graduate Texts in Mathematics},
    VOLUME = {21},
   EDITION = {Second},
 PUBLISHER = {Oxford University Press, Oxford},
      YEAR = {2011},
     PAGES = {xiv+684},
      ISBN = {978-0-19-960339-8},
   MRCLASS = {05-01 (05B35 90C27)},
  MRNUMBER = {2849819},
MRREVIEWER = {Maruti M. Shikare},
       DOI = {10.1093/acprof:oso/9780198566946.001.0001},
       URL = {https://doi.org/10.1093/acprof:oso/9780198566946.001.0001},
}

\end{document}